\documentclass[english]{amsart}

\usepackage{babel}
\usepackage{amstext}
\usepackage{amsmath}
\usepackage{amsfonts}
\usepackage{latexsym}
\usepackage{ifthen}
\usepackage{xypic}
\xyoption{all}
\newcommand{\id}{{\rm id}}
\renewcommand{\O}{{\mathcal O}}
\newcommand{\E}{{\mathcal E}}
\newcommand{\PN}{{\mathbb P}}
\newcommand{\rk}{{\rm rk \,}}
\newcommand{\Pic}{{\rm Pic}}
\newcommand{\lra}{\longrightarrow}
\newcommand{\KC}{{\mathbb C}}
\newcommand{\KZ}{{\mathbb Z}}

\newcommand{\End}{{\mathcal E}nd}
\newcommand{\D}{\mathbb D}

\newcounter{lemma}

\newtheorem{lemma1}[lemma]{\setcounter{equation}{0}}

\newenvironment{lemma}{\begin{lemma1}{\bf Lemma.}}{\end{lemma1}}

\newenvironment{theorem2}[1]{\begin{lemma1}{\bf Theorem [#1].}}{\end{lemma1}}
\newenvironment{proposition}{\begin{lemma1}{\bf Proposition.}}{\end{lemma1}}
\newenvironment{corollary}{\begin{lemma1}{\bf Corollary.}}{\end{lemma1}}
\newenvironment{remark}{\begin{lemma1}{\bf Remark.}\rm}{\end{lemma1}}
\newenvironment{definition}{\begin{lemma1}{\bf Definition.}}{\end{lemma1}}

\begin{document}

\title {Projective Threefolds with Holomorphic Conformal Structure} 
\author{Priska Jahnke}
\author{Ivo Radloff}
\address{Mathematisches Institut \\ Universit\"at Bayreuth \\ D-95440 Bayreuth/Germany}
\email{priska.jahnke@uni-bayreuth.de}
\email{ivo.radloff@uni-bayreuth.de}
\thanks{The authors gratefully acknowledge support by the Schwerpunkt program {\em Globale Methoden in der komplexen Geometrie} of the Deutsche Forschungsgemeinschaft.}
\date{\today}
\maketitle

%%%%%%%%%%%%%%%%%%%%%%%%%%%%

\section*{Introduction}
\setcounter{lemma}{0}

Complex manifolds modeled after projective space or hyperquadrics have been studied in connection with uniformization, twistor theory, and Fano geometry. Infinitesimal are the notions of holomorphic projective connections and holomorphic conformal structures (see \S 1-2 for the precise definitions). The complete list of surfaces with a holomorphic projective connection or a holomorphic conformal structure is due to Kobayashi and Ochiai. Projective threefolds with a holomorphic projective connection were completely classified by the authors in \cite{JR}. In this article we prove:

\vspace{0.2cm}

\noindent {\bf Theorem. }{\em The list of projective threefolds with a holomorphic conformal structure is as follows:
 \begin{enumerate}
   \item $Q_3$;
   \item {\'e}tale quotients of abelian threefolds;
   \item threefolds with universal covering space the threedimensional Lie ball $\D^{IV}_3$.
 \end{enumerate}}

\vspace{0.2cm}

\noindent We denote here by $\D^{IV}_3$ the bounded symmetric domain dual to the threedimensional hyperquadric $Q_3$. The above three types of manifolds are K\"ahler-Einstein; they all carry a flat holomorphic conformal structure in a canonical way (see \S 1). As in the surface case we have

\vspace{0.2cm}

\noindent {\bf Corollary~1. }{\em In projective dimension three, any holomorphic conformal structure is flat, i.e., comes from a quadric structure.}

\vspace{0.2cm}

A classical result identifies the threedimensional hyperquadric $Q_3$ with the Lagrangian Grassmannian parametrizing isotropic $2$-planes of a non degenerate skew symmetric bilinear form of rank four. As a consequence, $T_{Q_3}$ is the symmetric square of a rank two vector bundle, and we get

\vspace{0.2cm}

\noindent {\bf Corollary~2. }{\em The tangent bundle of a complex projective threefold $M$, or of some finite {\'e}tale covering of $M$, is a twisted symmetric square $T_M \simeq S^2\E \otimes {\mathcal L}$ if and only if $M$ is in the above list.}

\vspace{0.2cm}

\noindent Here of course $\E$ is some rank two vector bundle and ${\mathcal L} \in \Pic(M)$.

In the proof of the Theorem, Mori theory and elementary results from VHS are our main tools. The restriction to the projective case is owed to the fact that Mori Theory is yet only sufficiently settled for projective manifolds. 

The main Theorem together with \cite{JR} give a complete answer to the questions concerning manifolds modeled after hermitian symmetric spaces raised by Kobayashi and Ochiai in projective dimension three. 

\vspace{0.2cm}

\noindent {\bf Acknowledgements.} The authors are grateful to Professor Th. Peternell for valuable discussions and to Professor N. Mok for some help with Corollary~2.

%%%%%%%%%%%%%%%%%%%%

\section{Preliminaries}
\setcounter{lemma}{0}

We work over the complex numbers $\KC$. A manifold always is a compact connected holomorphic complex manifold. We will usually assume $M$ is K\"ahler.

\subsection{Manifolds modeled after Hyperquadrics} 
(cf. \cite{KobWu}) Let $Q_n$ be the hyperquadric in $\PN_{n+1}(\KC)$. In the sequel, we will usually think of $Q_n$ as the zero set of
  \[-2z_0z_{n+1} + \sum\nolimits_{i=1}^n z_i^2.\]
Think of $O(n+2)$ as the group of biholomorphic automorphisms of $Q_n$. A {\em quadric structure} is defined as follows.

Let $M$ be a complex manifold of dimension $n$ and $\{U_{\alpha}; z_{\alpha}^1, \dots, z_{\alpha}^n\}_{\alpha \in I}$ a holomorphic atlas of $M$. $M$ is said to have a {\em quadric structure}, if each $U_{\alpha}$ can be embedded into $Q_n$, such that the transition function is given by the restriction of a map in $O(n+2)$. Examples of manifolds with a quadric structure are
\begin{enumerate}
  \item $Q_n$;
  \item {\'e}tale quotients of abelian varieties; 
  \item $n$-folds with universal covering space $\D^{IV}_n$, the $n$-dimensional Lie ball.
\end{enumerate}
Recall that the $n$-dimensional Lie ball is the noncompact dual of the irreducible hermitian symmetric space $Q_n$, $(n \ge 3)$ (\cite{Helgason}). In adopted coordinates
\[\D^{IV}_n = \{z \in \KC^n \mid 2\sum |z_{\nu}|^2 < 1 + \big|\sum z_{\nu}^2\big|^2 < 2\},\]
and the embedding into $Q_n$ is given by 
  \[z \mapsto [1:z:\frac{1}{2}\sum z_i^2].\]
The above list of three examples will be refered to as {\em the list of standard examples}.

\

Kobayashi and Ochiai (following Chern, cf. \cite{KobNom}) introduced the concept of manifolds modeled after arbitrary hermitian symmetric spaces. There always is a corresponding infinitesimal notion.

%%%

\subsection{Holomorphic Conformal Structures}(cf. \cite{KobWu})
The infinitesimal notion to quadric structure is that of a holomorphic conformal structure. A holomorphic conformal structure is a special principal subbundle of the holomorphic frame bundle. The following definition is equivalent.

\begin{definition}
A manifold $M$ is said to have a holomorphic conformal structure, if there exists a nondegenerate holomorphic form
\[\omega \in H^0(M, S^2\Omega^1_M \otimes L),\]
where $L \in \Pic(M)$ is some line bundle.
\end{definition}

If $M$ carries a quadric structure, then it carries a holomorphic conformal structure. On $Q_n$, in homogeneous coordinates as above, the form is given by 
 \[-2dz_0dz_{n+1} + (dz_1)^2 + \cdots + (dz_n)^2\in H^0(Q_n, S^2\Omega_{Q_n}^1 \otimes \O_{\PN_{n+1}}(2)).\]
The form is $G \simeq O(n+2)$ invariant, pulling it back chart-wise we see that a manifold with a quadric structure indeed carries a holomorphic conformal structure. Holomorphic conformal structures on manifolds carrying a quadric structure are called {\em flat}.

\

From now on assume $M$ carries a holomorphic conformal structure, $n = \dim M$. Since $\omega$ is non degenerate, it induces an isomorphism
 \[T_M \simeq \Omega_M^1 \otimes L.\]
The determinant gives 
 \begin{equation} \label{detform}
   (\det T_M)^{\otimes 2} \simeq L^{\otimes n},
 \end{equation}
hence $\O_M(K_M)$ is divisible by $\frac{n}{2}$ in $\Pic(M)$. Certain relations hold among the Chern classes of $M$, coming from those of the quadric. For example, if $M$ is K\"ahler of dimension three, then
 \begin{equation} \label{chern}
   9c_2 = 4c_1^2 , \quad 27c_3 = 2c_1^3,
 \end{equation}
showing that the higher Chern classes may be computed from the first. For the general formulas see \cite{KobWu}.

\subsection{Known Results} 
The main results are due to Kobayashi and Ochiai, Ye, and Hwang and Mok. Kobayashi and Ochiai completely classified all compact complex surfaces, carrying a holomorphic conformal structure, in \cite{KoOcHQ}: 

\begin{theorem2}{Kobayashi Ochiai}
The list of compact complex surfaces (not only K\"ahler) admitting a holomorphic conformal structure is as follows:
\begin{enumerate}
  \item $Q_2 = \PN_1 \times \PN_1$;
  \item a ruled surface over a curve $C$ of genus $\ge 1$ such that the covering transformations of $\tilde{C} \times \PN^1$ respect the natural quadric structure on $\tilde{C} \times \PN^1$;
 \item a torus, a bielliptic surface or an elliptic surface with even $b_1$ and $c_1^2 = 0 = c_2$;
 \item a surface covered by the bidisc;
 \item an Inoue surface $S_U$ (see \cite{Inoue});
 \item a Hopf surface $\KC^{2*}/\Gamma$, where $\Gamma$ contains only transformations of the form $(z_1, z_2) \mapsto (\alpha z_1, \beta z_2)$ or $(z_1, z_2) \mapsto (\alpha z_2, \beta z_1)$.
\end{enumerate}
Any of these surfaces does already admit a quadric structure.
\end{theorem2}
A surface carries a holomorphic conformal structure if and only if the tangent bundle, after some finite {\'e}tale cover, splits as a sum of two line bundles. Indeed, $\omega$ induces a multiple section of
  \[\PN(\Omega_M^1) \lra M.\]
Depending on whether this section is irreducible or not, the tangent bundle of $M$ itself or of an {\'e}tale 2:1-cover of $M$ splits (alternatively notice $SO(2) \simeq \KC^*$). In arbitrary dimensions, the following is known:

\begin{theorem2}{Kobayashi-Ochiai}
 The list of K\"ahler-Einstein manifolds with a holomorphic conformal structure is the list of standard examples.
\end{theorem2}

\begin{theorem2}{Ye} \label{YeTheo}
If $M$ is a projective manifold with a holomorphic conformal structure of dimension $n \ge 3$ and $K_M$ is not nef, then $M$ is the hyperquadric. 
\end{theorem2}
In the case of a Fano manifold of odd dimension, this is just Kobayashi and Ochiai's criterion of hyperquadrics. For uniruled manifolds see also \cite{HM}.

%%%%%

\section{Holomorphic Conformal Strutures and Atiyah Classes}
\setcounter{lemma}{0}

We first recall the basic definition (cf. \cite{Ati}). 

\subsection{Atiyah Classes} 
Let $E$ be a holomorphic vector bundle of rank $r$ on $M$ and $\{U_{\alpha}; z_{\alpha}^1, \dots, z_{\alpha}^n\}_{\alpha \in I}$ holomorphic coordinate charts
 \footnote{After refinement, we will assume that $U_{\alpha_1} \cap \dots \cap U_{\alpha_k}$ are all Stein without further mentioning it.}
 on $M$ with coordinates $z_{\alpha}^i$, such that $E$ is trivial on $U_{\alpha}$. Let $\{U_{\alpha}; e_{\alpha}^1, \dots, e_{\alpha}^r\}_{\alpha \in I}$ be a local frame for $E$ and denote by $g_{\alpha\beta} \in H^0(U_{\alpha} \cap U_{\beta}, {\rm Gl}(r, \O_M))$ the corresponding transition functions of $E$ such that $e_{\beta}^k = \sum_l g_{\alpha\beta}^{lk}e_{\alpha}^l$.

The {\em Atiyah class of $E$} is the splitting obstruction of the first jet sequence
  \[0 \lra \Omega^1_M \otimes E \lra J_1(E) \lra E \lra 0,\]
i.e., it is the image of $\id_E$ under the first connecting morphism 
 \[H^0(M, \End(E)) \lra {\rm Ext}^1(E, \Omega_M^1 \otimes E) \simeq H^1(M, \End(E) \otimes \Omega^1_M).\] 
The Dolbeault isomorphism $H^1(M, \End(E) \otimes \Omega_M^1) \simeq H^{1,1}(M, \End(E))$ maps the Atiyah class to $[-\Theta_h]$, where $\Theta_h$ denotes the canonical curvature of $E$ with respect to a hermitian metric $h$ on $E$. In particular, the trace of the Atiyah class is $-2\pi i c_1(E)$ in $H^1(M, \Omega^1_M)$.

If we define $a(E)$ as $-\frac{1}{2\pi i}$ times the Atiyah class of $E$, then the trace of $a(E)$ is $c_1(E)$, which makes this definition convenient for our purposes. In local coordinates, $a(E)$ is the class of the \v{C}hech cocycle $a(E)_{\alpha\beta} \in {\mathcal C}^1({\mathcal U}, \End(E)\otimes \Omega_M^1)$, where
 \[a(E)_{\alpha\beta} = \frac{1}{2\pi i} \sum_{i,j,l} \frac{\partial g_{\alpha\beta}^{jl}}{\partial z_{\alpha}^i}dz_{\alpha}^i \otimes e_{\alpha}^j \otimes e_{\beta}^{l*} = \frac{1}{2\pi i} \sum_{1 \le j,l \le r} dg_{\alpha\beta}^{jl} e_{\alpha}^j \otimes e_{\beta}^{l*}.\]
See \cite{Ati} for the functorial behavior of $a(E)$ under pull-back, tensor products and direct sums.

\subsection{The Atiyah Class of $\Omega^1_M$} 
Now let $M$ again be a compact K\"ahler manifold carrying a holomorphic conformal structure. The form  $\omega \in H^0(M, S^2\Omega_M^1 \otimes L)$ induces naturally an isomorphism $T_M \otimes L^* \to \Omega_M^1$. Denote by 
  \[\omega^{-1}: \Omega_M^1 \stackrel{\sim}{\lra} T_M\otimes L^*\]
the invers map. Applying $\omega^{-1}$ to $c_1(K_M) \in H^1(M, \Omega^1_M)$ gives $\omega^{-1}(c_1(K_M)) \in H^1(M, T_M\otimes L^*)$. The tensor product of $\omega$ and $\omega^{-1}(c_1(K_M))$ gives a cocycle in $H^1(M, (S^2\Omega_M^1 \otimes L) \otimes (T_M \otimes L^*))$. Hence
  \[\omega \otimes \omega^{-1}(c_1(K_M)) \in H^1(M, S^2\Omega_M^1 \otimes T_M).\]
The Atiyah class of $\Omega_M^1$ is an element in $H^1(M, S^2\Omega_M^1 \otimes T_M)$. The bundle $S^2\Omega_M^1 \otimes T_M$ is a subbundle of $\Omega_M^1 \otimes T_M \otimes \Omega_M^1$ and we have the identities
  \[\End(\Omega_M^1) \otimes \Omega_M^1 \simeq \Omega^1_M \otimes T_M \otimes \Omega_M^1 \simeq \Omega_M^1 \otimes \End(\Omega_M^1).\]
Any class $\xi \in H^1(M, \Omega_M^1)$ induces two classes in $H^1(M, S^2\Omega_M^1 \otimes T_M)$, namely $\id_{\Omega_M^1} \otimes \xi$ and $\xi \otimes \id_{\Omega_M^1}$. We claim

\begin{proposition}
If $M$ is a compact K\"ahler manifold of dimension $n$ with a holomorphic conformal structure, then the identity
 \[a(\Omega_M^1) = \id_{\Omega_M^1} \otimes \frac{c_1(K_M)}{n} + \frac{c_1(K_M)}{n}\otimes \id_{\Omega_M^1} - \frac{1}{2in\pi}\omega \otimes \omega^{-1}(c_1(K_M))\]
holds in $H^1(M, S^2\Omega_M^1 \otimes T_M)$.
\end{proposition}
The formula is in fact just a reformulation of standard formulas from conformal geometry (see \cite{KobWu}, 4.3.). 

\begin{proof} Viewing $\omega$ as a tensor, but omitting the tensor symbol, we may write on every holomorphic chart $U_{\alpha}$
\[\omega = \sum \varpi_{\alpha, jk}dz_{\alpha}^jdz_{\alpha}^k \quad \mbox{and} \quad \omega^{-1} = \sum \varpi_{\alpha}^{ih}\frac{\partial}{\partial z_{\alpha}^i}\frac{\partial}{\partial z_{\alpha}^h}.\]
On $U_{\alpha\beta}$ we have $(\varpi_{\alpha, jk}) = f_{\alpha\beta}(\varpi_{\beta, jk})$, where $f_{\alpha\beta} \in \O^*_{U_{\alpha\beta}}$, representing the line bundle $L$. By (\ref{detform}) we may assume $f_{\alpha\beta}^n = (\det\frac{\partial z_{\alpha}}{\partial z_{\beta}})^2$. 

Write $\sigma_{\alpha\beta k} = \frac{\partial}{\partial z_{\beta}^k} (\log \det (\frac{\partial z_{\alpha}^i}{\partial z_{\beta}^j}))$, such that $\sum \sigma_{\alpha\beta k}dz_{\beta}^k$ represents the first Chern class of $K_M$ on $U_{\alpha \beta}$, and $\omega^{-1}(c_1(K_M)) \otimes \omega$ is represented by
\[\sum \varpi_{\beta}^{ih}\varpi_{\beta, jk}\sigma_{\alpha\beta h}\frac{\partial}{\partial z_{\beta}^i}dz_{\beta}^jdz_{\beta}^k.\]
The claim of the Proposition, the formula multiplied by $2\pi i$, is that there exists a solution to the cocycle $\xi \in {\mathcal C}^1(M, S^2\Omega_M^1 \otimes T_M)$ defined on $U_{\alpha\beta}$ as
  \[\xi_{\beta\alpha} = \sum\nolimits \left(\frac{\partial z_{\beta}^i}{\partial z_{\alpha}^a}\frac{\partial^2 z_{\alpha}^a}{\partial z_{\beta}^j \partial z_{\beta}^k} - \frac{1}{n}\big(\delta^i_j\sigma_{\alpha\beta k} + \delta^i_k\sigma_{\alpha\beta j} - \varpi_{\beta}^{ih}\varpi_{\beta jk} \sigma_{\alpha \beta h}\big)\right)\frac{\partial}{\partial z_{\beta}^i}dz_{\beta}^j dz_{\beta}^k.\]
A direct computation shows that a solution is given by $\zeta \in {\mathcal C}^0(M, S^2\Omega_M^1 \otimes T_M)$, on $U_{\alpha}$ defined as
\[\zeta_{\alpha} = \sum \frac{1}{2}\varpi_{\alpha}^{ih}\left(\frac{\partial \varpi_{\alpha, hj}}{\partial z_{\alpha}^k}  + \frac{\partial \varpi_{\alpha, hk}}{\partial z_{\alpha}^j} + \frac{\partial \varpi_{\alpha, jk}}{\partial z_{\alpha}^h}\right)\frac{\partial}{\partial z_{\alpha}^i}dz_{\alpha}^j dz_{\alpha}^k,\]
i.e., the coboundary of $\zeta$ is $\xi$. Note that both tensors are indeed symmetric.
\end{proof}

%%%%

\section{Abundant Manifolds}
\setcounter{lemma}{0}

Ye's result reduces the study of projective manifolds carrying a holomorphic conformal structure to those where $K_M$ is nef. The manifold is called minimal in this case. It is called {\em abundant} (or rather $K_M$ is called abundant), if $|mK_M|$ is spanned for some $m \gg 0$. The induced map  \[f: M \lra Y\]
after Stein factorization is called the {\em Iitaka fibration}.
If $M$ is abundant, then it is minimal. It is a conjecture that the convers holds as well. The conjecture is known to be true for $\dim M \le 3$. 

\begin{proposition}
 Let $M$ be a projective (or K\"ahler) abundant manifold with a holomorphic conformal structure. Let $Z$ be a smooth K\"ahler manifold, admitting a map $g: Z \to M$ such that $f \circ g$ maps $Z$ to a point. Then
  \[a(g^*\Omega_M^1) = 0 \; \mbox{ in } \; H^1(Z, \End(g^*\Omega_M^1)\otimes\Omega_Z^1).\]
\end{proposition}
  
\begin{proof}
Let $g: Z \to M$ be any holomorphic map. The Atiyah class of $g^*\Omega_M^1$ is obtained by first pulling back $a(\Omega_M^1) \in H^1(M, \End(\Omega_M^1)\otimes \Omega_M^1)$ to a class in $H^1(Z, \End(g^*\Omega_M^1)\otimes g^*\Omega_M^1)$, before applying the differential $dg$ to the last factor. Using the description of $\Omega_M^1$, we obtain
 \begin{equation} \label{at}
  a(g^*\Omega_M^1) = \frac{1}{n} \big(\id_{g^*\Omega_M^1} \otimes c_1(g^*K_M) + g^*c_1(K_M) \otimes dg + dg(\omega) \otimes g^*\omega^{-1}(g^*c_1(K_M))\big)
 \end{equation}
in $H^1(Z, \End(g^*\Omega_M^1)\otimes \Omega_Z^1)$. The last summand needs some explanation. We denote by $dg(\omega)$ the image of $g^*\omega$ under
\[dg: \;\; H^1(Z, g^*\Omega_M^1 \otimes g^*\Omega_M^1 \otimes g^*L) \lra H^1(Z, g^*\Omega_M^1 \otimes \Omega_Z^1 \otimes g^*L).\]
The second factor is obtained by viewing $g^*\omega^{-1}$ as a map $g^*\Omega_M^1 \to g^*T_M \otimes g^*L^*$, and by applying it to $g^*c_1(K_M) \in H^1(Z, g^*\Omega_M^1)$. The result is $g^*\omega^{-1}(g^*c_1(K_M)) \in H^1(Z, g^*T_M \otimes g^*L^*)$. In this way we get $dg(\omega) \otimes g^*\omega^{-1}(g^*c_1(K_M))$.

Now assume that $f \circ g$ maps $Z$ to a point. The map $f$ is given by $|mK_M|$, so $c_1(K_M)$ comes from $Y$. Then
  \[0 = g^*c_1(K_M) \in H^1(Z, g^*\Omega_M^1).\]
It is then clear that the first two summands in (\ref{at}) are zero. It is now also clear that $g^*\omega^{-1}$ applied to $g^*c_1(K_M)$ gives zero in $H^1(Z, g^*T_M \otimes g^*L^*)$. Then also the third summand vanishes and we conclude $a(g^*\Omega^1_M) = 0$.
\end{proof}

\begin{corollary} \label{ratcurves}
 None of the fibres of $f$ contains a rational curve.
\end{corollary}

\begin{proof}
 Assume that $g: \PN_1(\KC) \to M$ is a nontrivial map into some fibre of $f$. The proposition says $a(g^*\Omega_M^1) = 0$. Then $a(g^*T_M) = 0$ and since $\PN_1(\KC)$ is simply connected, $g^*T_M = \O_{\PN_1}^{\oplus n}$ (see \cite{Ati}). But $g^*T_M$ contains $T_{\PN_1} = \O_{\PN_1}(2)$. A contradiction.
\end{proof}

\begin{corollary} \label{equidim}
 The Iitaka fibration is equidimensional and the general fibre is covered by a torus.
\end{corollary}

\begin{proof}
 By a result of Kawamata, since $-K_M$ is $f$-nef, any fibre of $f$, whose dimension exceeds the dimension of the general fibre, is covered by rational curves. Corollary~\ref{ratcurves} implies that $f$ is equidimensional.

Let $F$ be the general fibre of $f$. We have the exact sequence
  \[0 \lra T_F \lra T_M|_F \lra N_{F/M} \lra 0.\]
The proposition gives $a(T_M|_F) = 0$, implying that all the Chern classes of $T_M|_F$ are zero. Since $N_{F/M} \simeq \O_F^{\oplus {\rm codim} F}$, all the Chern classes of $N_{F/M}$ are zero. Hence $c_i(T_F) = 0$ for all $i$. Then $F$ is covered by a torus.
\end{proof}

\begin{corollary} \label{big}
 If $K_M^n \not= 0$, then the universal cover of $M$ is $\D^{IV}_n$.
\end{corollary}

\begin{proof}
 If $K_M^n \not= 0$, then $f$ is a birational morphism. By Corollary~\ref{equidim}, $f$ is an embedding, implying that $K_M$ is ample. By the famous result of Yau, $M$ admits a K\"ahler-Einstein metric. By Kobayashi and Ochiai's result on K\"ahler-Einstein manifolds with holomorphic conformal struture, the universal cover of $M$ is $\D^{IV}_n$, the $n$-dimensional Lie ball.
\end{proof}

\begin{remark}
Abundant K\"ahler surfaces with holomorphic conformal structure in Kobayashi and Ochiai's list are tori, bielliptic surfaces, elliptic surfaces with $c_1^2 = c_2 = 0$ and surfaces covered by the bidisc. In all of these cases, the above results can be seen immediately.

In the case of a torus or bielliptic surface, $f$ maps to a point. For an elliptic surface, $f$ goes to a curve, and $c_2 = 0$ implies the only singular fibres are multiples of smooth fibres (\cite{Kodaira}). In the case of a surface covered by the bidisc, $f$ is an embedding.
\end{remark}

%%%%%%%%

\section{Proof of the Main Theorem -- Threefolds}
\setcounter{lemma}{0}

In order to prove the main Theorem, it suffices by Theorem~\ref{YeTheo} to study minimal projective threefolds $M$ with a holomorphic conformal structure. Then $M$ is abundant. Denote by
 \[f: M \lra Y\]
the Iitaka fibration as above. The dimension of $Y$ is given by the Kodaira dimension $\kappa(M) = \dim Y$. If $\kappa(M) = 0$, then $M$ is covered by an abelian variety according to Corollary~\ref{equidim}. If $\kappa(M) = 3$, then $K_M^3 \not= 0$ and $M$ is a quotient of $\D^{IV}_3$ by Corollary~\ref{big}. It remains to show $\kappa(M) \not= 1, 2$. The proof of our main Theorem is concluded in the following two sections:

%%%%

\subsection{Iitaka Fibration to a Surface}

\begin{proposition} \label{Surface}
 A (minimal) threefold of Kodaira dimension two does not carry a holomorphic conformal structure.
\end{proposition}

\begin{proof}
 Assume to the contrary that $M$ is a minimal threefold, carrying a holomorphic conformal structure, and $\kappa(M) = 2$. Then the Iitaka fibration
 \[f: M \lra S\]
is an equidimensional fibration onto a normal surface $S$. The general fibre of $f$ is a smooth elliptic curve. Corollary~\ref{ratcurves} implies that the $j$-invariant of the fibration is constant. The only singular fibres (of an elliptic surface swept out by the fibres over a general hyperplane section in $S$) are of type ${}_mI_0$.

First assume that $f$ is an elliptic bundle. Then after a finite \'etale cover we may assume that we have the exact sequence
 \begin{equation} \label{relseq}
  0 \lra f^*\Omega^1_S \lra \Omega_M^1 \lra \O_M \lra 0.
 \end{equation}
The corresponding long exact sequence of cohomology is exact on $H^0$-level, implying that (\ref{relseq}) is split exact. Since $L^{\otimes 3} \simeq \O_M(-2K_M)$ we have $H^0(M, L) = 0$. Then $\Omega_M^1 \otimes L \simeq T_M$ implies that $\O_M$ is a direct summand of $f^*\Omega^1_S \otimes L$. Then $c_2(f^*\Omega^1_S \otimes L) = 0$. In other words $c_2(f^*\Omega_S^1) = \frac{2}{9}c_1^2(M)$. Then the above sequence gives $c_2(\Omega^1_M) = c_2(f^*\Omega_S^1)$, i.e., $c_2(\Omega^1_M) = \frac{2}{9}c_1^2(M)$. Formula (\ref{chern}) implies $c_1^2(M) = 0$, contradicting $\kappa(M) = 2$.

In the general case, we still have a relative one form on $M$, and $\Omega_M^1$ has a trivial direct summand. The same computations show $c_2(\Omega_M^1) = \frac{2}{9}c_1^2(M)$ and therefore $c_1^2(M) = 0$, contradicting $\kappa(M) = 2$.
\end{proof}
 
%%%%

\subsection{Iitaka Fibration to a Curve}
Now assume that the Iitaka fibration $f: M \to C$ goes to a curve. Corollary~\ref{equidim} says the central case we have to deal with is a smooth abelian fibration. In the smooth case we have the exact sequence
  \[0 \lra f^*\Omega_C^1 \lra \Omega_M^1 \lra \Omega_{M/C}^1 \lra 0.\]
Put $E^{1, 0} = f_*\Omega^1_{M/C}$. This is a rank two vector bundle on $C$, and $f^*E^{1, 0} = \Omega^1_{M/C}$. Denote by $N^{1,0}$ the kernel of
  \[\theta^{1, 0}: E^{1, 0} = f_*\Omega^1_{M/C} \lra R^1 f_*f^*\Omega^1_C.\]
There are Arakelov type inequalities relating the degrees of these bundles. By \cite{ViZu}, $\deg N^{1, 0} \le 0$, which is all we need.

\begin{proposition} \label{Curve}
 A (minimal) threefold of Kodaira dimension one does not carry a holomorphic conformal structure.
\end{proposition}

\begin{proof}
 Assume to the contrary that $M$ is a minimal threefold with holomorphic conformal structure and $\kappa(M) = 1$. Corollary~\ref{equidim} implies that the general fibre $F$ is either abelian or hyperelliptic. None of the fibres of $f$ contains a rational curve. By \cite{Oguiso}, Theorem~B.1., the only singular fibres of $f$ are multiple fibres.

\vspace{0.2cm}

We first assume that $f$ is smooth. If $F$ is hyperelliptic, we have the so called intermediate Albanese, i.e., the Albanese map of $F$ gives rise to a commutative diagram
  \[M \lra Y' \lra C,\]
with $Y'$ a smooth surface. The map $f': M \to Y'$ is a smooth elliptic fibration, which is excluded by Proposition~\ref{Surface}. If $F = A$ is an abelian surface, then we conclude as follows. From
  \[0 \lra T_A \simeq \O_A^{\oplus 2} \lra T_M|_A \lra N_{A/M} \simeq \O_A \lra 0\]
we see $h^0(A, T_M|_A) \ge 2$, and $h^0(A, \Omega_M^1|_A) \ge 1$. Tensorize with $L^*$ and use $\Omega_M \simeq T_M \otimes L^*$ to get $H^0(A, L^*|_A) \not= 0$. Some multiple of $L$ is in $f^*\Pic(C)$, hence $L|_A$ is torsion. We conclude $L = f^*L_C$, where $L_C = f_*L \in \Pic(C)$. Then
  \[f_*\Omega^1_M \simeq f_*T_M \otimes L^*_C \quad \mbox{and} \quad 2 \le \rk f_*T_M = \rk f_*\Omega^1_M \le 3.\]

Assume $\rk f_*\Omega^1_M = 2$. Then $f_*T_M \simeq (E^{1, 0})^*$ and $f_*\Omega_M^1 \simeq (E^{1, 0})^* \otimes L_C^*$ and we get the sequence
  \[0 \lra \Omega_C^1 \lra (E^{1, 0} \otimes L_C)^* \simeq f_*\Omega_M^1 \lra N^{1, 0} \lra 0,\]
implying $\deg K_C = -\deg (E^{1, 0} \otimes L_C) - \deg N^{1,0}$. Using $K_M = K_{M/C} + f^*K_C$ and $L^{\otimes 3} \simeq \O_M(-2K_M)$ we conclude $\frac{1}{3}K_M \equiv f^*N^{1,0}$. Since $N^{1,0}$ is seminegative, this is a contradiction.

Hence $\rk f_*\Omega^1_M = 3$. Then $f_*\Omega_{M/C} = E^{1, 0} = N^{1, 0}$ is flat (\cite{Kollar}, Proposition~4.10) and $K_M \equiv f^*K_C$ (simplest example being $X = A \times C$). From $K_M \equiv f^*K_C$ we infer $C \not\simeq \PN_1$. From
  \[0 \lra \Omega_C^1 \lra f_*\Omega_M^1 \simeq f_*T_M \otimes L^*_C \lra E^{1, 0} \lra 0\]
we conclude that $f_*T_M \otimes L^*_C$ is nef. From the push-forward of the relative tangent sequence tensorized with $L_C^*$, we get a surjection $f_*T_M \otimes L^*_C \to T_C \otimes L_C^*$ and conclude $\deg (T_C \otimes L_C^*) \ge 0$. Using $K_M \equiv f^*K_C$ and $L^{\otimes 3} \simeq \O_M(-2K_M)$ once again we obtain $-\frac{1}{3}\deg K_C \ge 0$. Then $C$ must be elliptic, implying $K_M \equiv 0$, again a contradiction.

\vspace{0.2cm}

If $f$ is not smooth, then $f$ is an almost smooth fibration with multiple fibres $F_1, \dots, F_l$ of multiplicities $m_1, \dots, m_l$. Define
\[D = \sum\nolimits_{i=1}^l (m_i-1)F_i.\] 
A ramified base-change $C' \to C$ of degree $m = {\rm l.c.m.}(m_1, \dots, m_l)$, ramified over the critical values of $f$ and one additional point $a_0 \in C$, leads to a smooth fibration $f': M' \to C'$, where $M'$ is the normalization of $M \times_C C'$ (see for example \cite{Kodaira}). Let $\mu: M' \to M$ be the induced ramified cover. On $M'$ we have the commutative diagram
 \[\xymatrix{0 \ar[r] & \mu^*(f^*K_C \otimes \O_M(D)) \ar[r] \ar@{^{(}->}[d] & \mu^*\Omega_M^1 \ar[r]^{\mbox{$\rho$}} \ar@{^{(}->}[d] & \Omega_{M'/C'}^1 \ar[r] \ar@{=}[d] & 0\\
 0 \ar[r] & {f'}^*K_{C'} \ar[r] & \Omega_{M'}^1 \ar[r] & \Omega_{M'/C'}^1 \ar[r] & 0,}\]
where $\rho$ is the pull-back of the cokernel map of $f^*K_C \otimes \O_M(D) \hookrightarrow \Omega_M^1$. Note that $\Omega_{M'/C'}^1$ is the pull-back of the rank two vector bundle ${E^{1,0}}' = f'_*\Omega_{M'/C'}^1$ and $\mu$ is unramified over the general fibre of $f$. With $\mu^*T_M \simeq \mu^*\Omega_M^1 \otimes \mu^*L$, essentially the same argument as above shows the claim, using the push-forward of both of the above sequences. 
\end{proof}

%%%%%%%%%%%%%%%%%

\section{Symmetric Squares}
\setcounter{lemma}{0}

In the case of surfaces, the existence of a holomorphic conformal structure is equivalent to the splitting of the tangent bundle as a sum of two line bundles. In the threedimensional case, an alternative description is given in Corollary~2, which we are now going to prove:

\begin{proof}[Proof of Corollary~2] If $M$ is a projective threefold and $T_M \simeq S^2\E \otimes {\mathcal L}$, with $\E$ some rank two vector bundle and ${\mathcal L} \in \Pic(M)$, then the canonical isomorphism 
 \[\E \simeq \E^* \otimes \det\E\]
implies $T_M \simeq \Omega_M^1 \otimes L$, where $L = \det\E^{\otimes 2} \otimes {\mathcal L}^{\otimes 2}$, and this isomorphism induces a non degenerate $\omega \in H^0(M, S^2\Omega_M^1 \otimes L)$. Hence $M$ carries a holomorphic conformal structure. 

For the convers let $M$ be a projective threefold with a holomorphic conformal structure. We have to show that there exists some finite {\'e}tale covering $\tilde{M} \to M$ such that $T_{\tilde{M}}$ is a (twisted) symmetric square. With $F = T_M \otimes L \otimes (\det T_M)^*$ this asks for whether we can lift the class in $H^1(M, SO(3, \KC))$ defined by $F$ to $H^1(M, SL(2, \KC))$ using
  \[1 \lra \KZ_2 \lra SL(2, \KC) \lra SO(3, \KC) \lra 1.\]
The obstruction lies in $H^2(M, \KZ_2)$. However, using the main Theorem, we can conclude case by case. 

The case $M = Q_3$ is clear. The threedimensional quadric is a Lagrangian Grassmannian, hence $T_{Q_3} \simeq S^2U$ where $U$ has rank two. The case $M$ a quotient of an abelian threefold is trivial. The case where the universal covering space of $M$ is $\D^{IV}_3$ perhaps needs some explanations. Think of $\D^{IV}_3$ as embedded into $Q_3$ and restrict $U$ and $T_{Q_3}$ to $\D^{IV}_3$. As soon as $U$ descends to $M$, the isomorphism $T_{Q_3} \simeq S^2U$ descends, and $T_M$ is a symmetric square. Recall the sequence
  \[1 \lra \KZ_2 \lra Sp(2, \KC) \lra SO(5, \KC) \lra 1.\]
The bundle $U$ is not $SO(5, \KC)$, but $Sp(2, \KC)$ homogeneous; the problem is $-\id \in Sp(2, \KC)$ which acts trivially on $\D^{IV}_3$, but not on $U$. A sufficient condition for $U$ to descend to a vector bundle on $M$ therefore is that the group by which we divide $\D^{IV}_3$ does not contain $-\id$. 

The group is residual finite. If it contains $-\id$, we can choose a subgroup of finite index, which does not contain $-\id$. The quotient of $\D^{IV}_3$ by this group gives a finite {\'e}tale covering $\tilde{M} \to M$ such that $T_{\tilde{M}}$ is a symmetric square.
\end{proof}

%%%%%%%%%%%%%%%%%%

\end{document}